\documentclass[11pt]{amsart}
\usepackage{amsmath,graphicx}
\usepackage[all]{xy}

\newtheorem{thm}{Theorem}[section]
\newtheorem{lem}[thm]{Lemma}
\newtheorem{prop}[thm]{Proposition}
\newtheorem{cor}[thm]{Corollary}

\theoremstyle{definition}

\theoremstyle{remark}
\newtheorem{remark}[thm]{Remark}

\theoremstyle{plain}

\newcommand{\Z}{{\mathbb{Z}}}

\newcommand{\R}{{\mathbb{R}}}
\newcommand{\Zf}{\Z^2_{\mathrm{pr}}}
\newcommand{\N}{{\mathbb{N}}}

\newcommand{\FF}{{\mathcal{F}}}
\newcommand{\TT}{{\mathcal{T}}}
\newcommand{\DD}{{\mathcal{D}}}
\newcommand{\odd}{{\mathrm{odd}}}
\newcommand{\even}{{\mathrm{even}}}
\newcommand{\area}{{\mathrm{Area}}}
\newcommand{\lh}{{\mathrm{length}}}

\numberwithin{equation}{section}

\begin{document}

\title[The Farey sequence with odd denominators]
{On the distribution of the Farey sequence with odd denominators}

\author[F.P. Boca, C. Cobeli, A. Zaharescu]{Florin P. Boca, Cristian Cobeli
and Alexandru Zaharescu}

\address{FPB and AZ: Department of Mathematics, University of Illinois,
1409 W. Green Str., Urbana IL 61801, USA}

\address{FPB, CC and AZ: Institute of Mathematics of the
Romanian Academy, P.O. Box 1-764, RO-014700 Bucharest, Romania}

\address{fboca@math.uiuc.edu, cristian.cobeli@imar.ro,
zaharesc@math.uiuc.edu}

\thanks{Research partially supported by ANSTI grant C6189/2000}


\maketitle

\section{Introduction and Statement of Results}
Given a positive integer $Q$, we denote by $\FF_Q$ the set of
irreducible rational fractions in $(0,1]$ whose denominators do
not exceed $Q$. That is,
\begin{equation*}
\FF_Q =\{ a/q \, ;\, 1\leq a\leq q\leq Q,\ \gcd (a,q)=1\} .
\end{equation*}

Problems concerning the distribution of Farey fractions have been
studied in the 20's by Franel and Landau (\cite{Fra},\cite{La})
and more recently in \cite{ABCZ}, \cite{BCZ0}, \cite{BCZ},
\cite{BGZ}, \cite{Hall1}, \cite{Hall2}, \cite{HT}, \cite{HSZ},
\cite{Hay}, \cite{Hux}, \cite{HZ}.

It is well-known that
\begin{equation*}
N_Q =\# \FF_Q =6Q^2 /\pi^2 +O(Q\log Q).
\end{equation*}
We denote by $\FF_Q^<$ the set of pairs $(\gamma,\gamma^\prime)$
of consecutive elements in $\FF_Q$.

In this paper we are concerned with the set
\begin{equation*}
\FF_{Q,\odd} =\{ a/q \in \FF_Q \, ;\, q \ \odd \}
\end{equation*}
of Farey fractions of order $Q$ with odd denominators. For
instance,
\begin{equation*}
\begin{split}
& \FF_8=\bigg\{ \frac{1}{8},\frac{1}{7},\frac{1}{6}, \frac{1}{5},
\frac{1}{4},\frac{2}{7},\frac{1}{3}, \frac{3}{8},\frac{2}{5},
\frac{3}{7},\frac{1}{2},
\frac{4}{7},\frac{3}{5},\frac{5}{8},\frac{2}{3},
\frac{5}{7},\frac{3}{4}, \frac{4}{5},\frac{5}{6},
\frac{6}{7},\frac{7}{8},1\bigg\}  \quad \mbox{\rm and} \\
& \FF_{8,\odd} =\bigg\{ \frac{1}{7},\frac{1}{5},
\frac{2}{7},\frac{1}{3},\frac{2}{5},\frac{3}{7},
\frac{4}{7},\frac{3}{5},\frac{2}{3},
\frac{5}{7},\frac{4}{5},\frac{6}{7},1\bigg\} .
\end{split}
\end{equation*}

The set of pairs $(\gamma,\gamma^\prime)$ of consecutive elements
in $\FF_{Q,\odd}$ is denoted by $\FF_{Q,\odd}^<$. It is not hard
to prove (see \cite{Hay}) that
\begin{equation}\label{1.1}
N_{Q,\odd} =\# \FF_{Q,\odd}=2Q^2 /\pi^2+O(Q\log Q).
\end{equation}

It is well-known that $\Delta (\gamma,\gamma^\prime):= a^\prime
q-aq^\prime =1$ whenever $\gamma =\frac{a}{q} <
\frac{a^\prime}{q^\prime} =\gamma^\prime$ are consecutive elements
in $\FF_Q$. This certainly fails when $\gamma <\gamma^\prime$ are
consecutive in $\FF_{Q,\odd}$. A first step in the study of the
distribution of the values of $\Delta (\gamma,\gamma^\prime)$ for
pairs $(\gamma,\gamma^\prime)$ of consecutive fractions in
$\FF_{Q,\odd}$ was undertaken by A. Haynes in \cite{Hay}. He
proved that if one denotes
\begin{equation*}
N_{Q,\odd} (k)=\# \{ \gamma <\gamma^\prime \ \mbox{\rm succesive
in} \ \FF_{Q,\odd} \, ;\, \Delta (\gamma,\gamma^\prime)=k\} ,
\end{equation*}
then the asymptotic frequency
\begin{equation*}
\rho_\odd (k) =\lim\limits_{Q\rightarrow \infty} \frac{N_{Q,\odd}
(k)}{N_{Q,\odd}}
\end{equation*}
exists, and is expressed as
\begin{equation*}
\rho_\odd (k)=\frac{4}{k(k+1)(k+2)}\, , \qquad k\in \N^* .
\end{equation*}

This can be written as
\begin{equation}\label{1.2}
\rho_\odd (k)=\begin{cases}
\area (\TT_k) & \mbox{\rm if $k\geq 2$} \\
\frac{1}{2}+\area (\TT_1) & \mbox{\rm if $k=1$,}
\end{cases}
\end{equation}
where (as in \cite{BCZ}) we denote $\TT_k =\big\{ (x,y)\in \TT \,
;\, \big[ \frac{1+x}{y} \big] =k\big\}$, $k\in \N^*$, and $\TT=\{
(x,y)\in [0,1] \, ;\, x+y>1\}$.

In this note we study, for fixed $h\geq 1$, the distribution of
consecutive elements $\gamma_i <\gamma_{i+1} <\cdots
<\gamma_{i+h}$ in $\FF_{Q,\odd}$, and compute the probability that
such an $(h+1)$-tuple satisfies $\Delta
(\gamma_i,\gamma_{i+1})=\Delta_1,\dots, \Delta
(\gamma_{i+h-1},\gamma_{i+h})=\Delta_h$. More precisely, we prove
that if one denotes
\begin{equation*}
N_{Q,\odd} (\Delta_1,\dots,\Delta_h)=\# \left\{
\begin{matrix} i ; & \hspace{-6pt} \gamma_i <\gamma_{i+1} <
\cdots <\gamma_{i+h} \ \mbox{\rm consecutive in} \ \FF_{Q,\odd} \\
& \Delta (\gamma_{i+j-1},\gamma_{i+j})= \Delta_j,\ j=1,\dots,h
\end{matrix} \right\} ,
\end{equation*}
then
\begin{equation*}
\rho_\odd (\Delta_1,\dots ,\Delta_h) =\lim\limits_{Q\rightarrow
\infty} \frac{N_{Q,\odd} (\Delta_1,\dots,\Delta_h)}{N_{Q,\odd}}
\end{equation*}
exists for all $h\geq 2$, and give an explicit formula for it.

To state the main result, we shall employ the area-preserving
transformation $T$ of $\TT$, introduced in \cite{BCZ} and defined
by
\begin{equation}\label{I3}
T(x,y)=\bigg( y,\bigg[ \frac{1+x}{y} \bigg] y-x\bigg).
\end{equation}
We denote
\begin{equation*}
\TT_{k_1,\dots,k_h}=\TT_{k_1} \cap T^{-1} \TT_{k_2} \cap \cdots
\cap T^{-h+1} \TT_{k_h} .
\end{equation*}

We notice that if $\gamma =\frac{a}{q} <\gamma^\prime
=\frac{a^\prime}{q^\prime} <\gamma^{\prime \prime} =
\frac{a^{\prime \prime}}{q^{\prime \prime}}$ are consecutive
elements in $\FF_Q$, then $T\big( \frac{q}{Q}, \frac{q^\prime}{Q}
\big) =\big( \frac{q^\prime}{Q}, \frac{q^{\prime \prime}}{Q}
\big)$. Moreover, if we set $\kappa (x,y)=\big[ \frac{1+x}{y}
\big]$, then the positive integer $\kappa \big(
\frac{q}{Q},\frac{q^\prime}{Q} \big) =\big[ \frac{Q+q}{q^\prime}
\big]$ coincides with the index $\nu_Q (\gamma)$ of the Farey
fraction $\gamma$ in $\FF_Q$ considered in \cite{HSZ}.

\begin{figure}[ht]
\begin{center}
$$\mbox{\small $\displaystyle \xymatrix@M=1pt{
 &   &  &   &   &
\star=O \ar@{-} [dl]_(0.8){\Delta_1=k_1}|{k_1} \ar@{-}
[ddrrr]^(0.8){\Delta_1=1}|{k_1}
&  &  &  & \\
 &   &  &   & E \ar@{-} [dl]|{k_2} &  &   &  &  & \\
 & & & O \ar@{-} [dl]_(0.8){k_2\, even}|{k_3}
\ar@{-} [dd]^(0.3){\substack{k_2 \, odd \\ \Delta_2=1}}|{k_3} & &
& & & O \ar@{-} [dl]_(0.8){\substack{k_1\, odd \\
\Delta_2=k_2}}|{k_2}
\ar@{-} [dd]^(0.3){\substack{k_1 \, even \\ \Delta_2=1}}|{k_2}  & \\
 &   & E \ar@{-} [dl]_(0.8){\Delta_2=k_3}|{k_4} & & & & &
E \ar@{-} [dl]|{k_3} & &  \\
 & O \ar@{-} [dl]_(0.4){k_4 \, even}|{k_5}
\ar@{-} [dd]^(0.25){\substack{k_4 \, odd\\ \Delta_3=1}}|{k_5} &  &
O \ar@{-} [dl]_(0.4){k_3\, odd}|{k_4} \ar@{-}
[dd]^(0.3){\substack{k_3 \, even \\ \Delta_3=1}}|{k_4} &  &  & O
\ar@{-} [dl]_(0.4){\substack{k_3\, even \\ \Delta_3=k_4}}|{k_4}
\ar@{-} [dd]^(0.3){\substack{k_3\, odd\\ \Delta_3=1}}|{k_4} &  & O
\ar@{-} [dl]_(0.4){\substack{k_2\, odd \\ \Delta_3=k_3}}|{k_3}
\ar@{-} [dd]^(0.3){\substack{k_2 \, even \\ \Delta_3=1}}|{k_3}  & \\
E \ar@{-} [d]_(0.25){\Delta_3=k_5}|{k_6} & & E \ar@{-}
[d]_(0.25){\Delta_3=k_4}|{k_5} & & &
E \ar@{-} [d]|{k_5} & & E \ar@{-} [d]|{k_4} & & \\
O & O & O & O & & O & O & O & O & }$}$$
\end{center}
\caption{The tree ${\mathfrak{T}}_3$} \label{Figure1}
\end{figure}
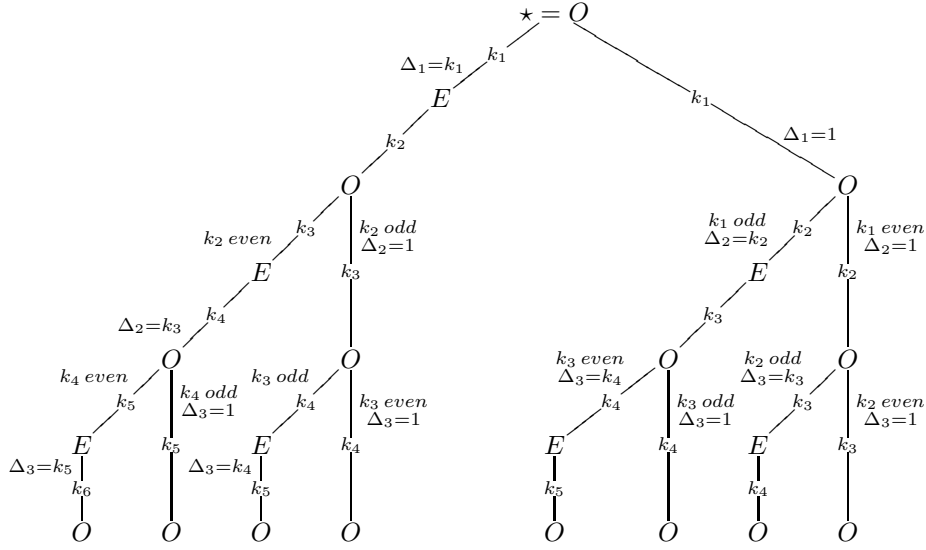

It will be worthwhile to consider the tree ${\mathfrak{T}}_h$
defined by the following properties:
\begin{itemize}
\item vertices are labeled by $O$ and $E$; \item the starting
vertex $\star$ is labeled by $O$; \item there is exactly one edge
starting from an $E$ vertex, and such an edge always ends into an
$O$ vertex; \item there are exactly two edges starting from an $O$
vertex, and they end (respectively) into an $E$ vertex and into an
$O$ vertex; \item the number of $O$ vertices (besides $\star$) on
any path that originates at $\star$ is equal to $h$.
\end{itemize}
See Figure \ref{Figure1}.

We also consider the set ${\mathfrak{L}}_h$ of labeled paths
\begin{equation*}
w=\Big( \xymatrix{\star=O \ar@{-} [r]^(0.6){k_1} & v_1 \ar@{-}
[r]^(0.55){k_2} & v_2 \ar@{-} [r]^(0.5){k_3} & \cdots \ar@{-}
[r]^(0.55){k_{\vert w\vert}} & v_{\vert w\vert} }\Big),\qquad k_j
\in \N^* .
\end{equation*}
on the tree ${\mathfrak{T}}_h$ that start at $\star$ and pass
through $h+1$ vertices labeled by $O$ (including $\star$). That
is, $\# \{ j\, ;\, v_j=O\}=h$. We set $o(O)=\odd$, $o(E)=\even$.

For each labeled path $w\in {\mathfrak{L}}_h$ and each $h$-tuple
$\Delta=(\Delta_1,\dots,\Delta_h)\in (\N^*)^h$ we define
$c_{OE}(w)$ and $c_\Delta (w)$ by induction as follows:
\begin{equation*}
\begin{split}
& c_{OE} \big( \xymatrix{\star=O \ar@{-} [r]^(0.6){k_1} & E
\ar@{-} [r]^(0.5){k_2} & O } \big)=k_1 , \qquad c_{OE} \big(
\xymatrix{\star=O \ar@{-} [r]^(0.6){k_1} & O } \big)=\emptyset ,\\
& c_{\Delta_1} \big( \xymatrix{\star=O \ar@{-} [r]^(0.6){k_1} & E
\ar@{-} [r]^(0.5){k_2} & O } \big) =\Delta_1 ,\qquad c_{\Delta_1}
\big( \xymatrix{\star=O \ar@{-} [r]^(0.6){k_1} & O
}\big)=\emptyset ,
\end{split}
\end{equation*}
For $w=w^\prime w^{\prime \prime} \in {\mathfrak{L}}_{h+1}$ with
$w^\prime \in {\mathfrak{L}}_h$ and $w^{\prime \prime}=\xymatrix{O
\ar@{-} [r]^(0.5){k} & E \ar@{-} [r]^(0.5){l} & O}$ or $w^{\prime
\prime}=\xymatrix{O \ar@{-} [r]^(0.5){k} & O}$, we have
\begin{equation*}
\begin{split}
& c_{OE}(w)=\begin{cases} (c_{OE}(w^\prime),k) & \mbox{\rm if
$w^{\prime \prime}=\xymatrix{O \ar@{-} [r]^(0.5){k} &
E \ar@{-} [r]^(0.5){l} & O}$} \\
c_{OE}(w^\prime) & \mbox{\rm if $w^{\prime \prime}=
\xymatrix{O \ar@{-} [r]^(0.5){k} & O}$}, \end{cases} \\
& c_{(\Delta_1,\dots,\Delta_{h+1})} (w) =\begin{cases} (c_\Delta
(w^\prime),\Delta_{h+1}) & \mbox{\rm if $w^{\prime \prime}=
\xymatrix{O \ar@{-} [r]^(0.5){k} &
E \ar@{-} [r]^(0.5){l} & O}$} \\
c_\Delta (w^\prime) & \mbox{\rm if $w^{\prime \prime}=\xymatrix{O
\ar@{-} [r]^(0.5){k} & O}$}. \end{cases}
\end{split}
\end{equation*}
For instance, if $w$ is the labeled path
\begin{equation*}
\xymatrix{\star=O \ar@{-} [r]^(0.6){k_1} & E \ar@{-}
[r]^(0.6){k_2} & O \ar@{-} [r]^(0.6){k_3} & O \ar@{-}
[r]^(0.6){k_4} & E \ar@{-} [r]^(0.6){k_5} & O \ar@{-}
[r]^(0.6){k_6} & O \ar@{-} [r]^(0.6){k_7} & E \ar@{-}
[r]^(0.6){k_8} & O}
\end{equation*}
in ${\mathfrak{L}}_5$, then
\begin{equation*}
c_{OE}(w)=(k_1,k_4,k_7) \qquad \mbox{\rm and} \qquad
c_{(\Delta_1,\dots,\Delta_5)} (w)=(\Delta_1,\Delta_3,\Delta_5).
\end{equation*}

We also denote by ${\mathfrak{S}}_\Delta$ the set of labeled paths
$$\xymatrix{v_0=\star =O \ar@{-} [r]^(0.7){k_1} &
v_1 \ar@{-} [r]^(0.5){k_2} & \cdots  \ar@{-} [r]^(0.5){k_{\vert
w\vert}} & v_{\vert w\vert}}$$
such that $c_{OE}(w)=c_\Delta (w)$,
and such that $k_j$ is even whenever it occurs as $\xymatrix{ E
\ar@{-} [r]^(0.5){k_j} & O \ar@{-} [r] & E}$ or as $\xymatrix{O
\ar@{-} [r]^(0.5){k_j} & O \ar@{-} [r] & O}$, and respectively odd
whenever it occurs as $\xymatrix{E \ar@{-} [r]^(0.5){k_j} & O
\ar@{-} [r] & O}$ or as $\xymatrix{O \ar@{-} [r]^(0.5){k_j} & O
\ar@{-} [r] & E}$.

Having established this notation, we may state our main result.

\medskip

\begin{thm}\label{MainThm}
Let $h\geq 1$, and let $\Delta=(\Delta_1,\dots ,\Delta_h) \in
(\N^*)^h$. Then
\begin{equation*}
\rho_{Q,\odd} (\Delta):= \frac{N_{Q,\odd}
(\Delta_1,\dots,\Delta_h)}{N_{Q,\odd}} =\rho_{\odd}(\Delta)+O_h
\bigg( \frac{\log^2 Q}{Q} \bigg)
\end{equation*}
as $Q\rightarrow \infty$, where
\begin{equation}\label{I4}
\rho_\odd (\Delta) =\sum\limits_{w\in {\mathfrak{L}}_h \cap
{\mathfrak{S}}_\Delta} \hspace{-5pt} \area (\TT_{k_1,\dots
,k_{\vert w\vert -1}}).
\end{equation}
\end{thm}

For $h=1$, this gives
\begin{equation*}
\rho_\odd (\Delta_1)=\begin{cases} \sum\limits_{k_1} \area
(\TT_{k_1})+ \area (\TT_1) =\frac{1}{2}+\area (\TT_1) &
\mbox{\rm if $\Delta_1 =1$} \\
\sum\limits_{k_2} \area (\TT_{\Delta_1} \cap T^{-1} \TT_{k_2})=
\area (\TT_{\Delta_1}) & \mbox{\rm if $\Delta_1 \geq 2$,}
\end{cases}
\end{equation*}
that is, the aforementioned result of Haynes.

For $h=2$, we obtain the following
\medskip

\begin{cor}\label{h2}
$\rho_{Q,\odd} (\Delta_1,\Delta_2)$ tends to $\rho_\odd
(\Delta_1,\Delta_2)$ for any $\Delta_1,\Delta_2 \in \N^*$ as
$Q\rightarrow \infty$. Moreover, we have:

{\em (i)} $\displaystyle \rho_\odd (1,1) =\sum\limits_{k_1 \,
\even} \area (\TT_{k_1}) +\sum\limits_{k_1 \, \odd} \area
(\TT_{k_1,1} ) +\sum\limits_{k_2 \, \odd} \area (\TT_{1,k_2}) $

$\displaystyle \qquad \qquad \qquad \qquad +\sum\limits_{k_2 \,
\even} \area (\TT_{1,k_2,1} ).$

{\em (ii)} if $\Delta_2 \geq 2$, then
\begin{equation*}
\rho_\odd (1,\Delta_2) =\sum\limits_{k_1 \, \odd} \area
(\TT_{k_1,\Delta_2} )+ \sum\limits_{k_2 \, \even} \area
(\TT_{1,k_2,\Delta_2}) .
\end{equation*}

{\em (iii)} if $\Delta_1 \geq 2$, then
\begin{equation*}
\rho_\odd (\Delta_1,1)=\sum\limits_{k_2 \, \odd} \area
(\TT_{\Delta_1,k_2} )+ \sum\limits_{k_2 \, \even} \area
(\TT_{\Delta_1,k_2,1} ).
\end{equation*}

{\em (iv)} if $\min (\Delta_1,\Delta_2)\geq 2$, then
\begin{equation*}
\rho_\odd (\Delta_1,\Delta_2)=\sum\limits_{k_2 \, \even} \area
(\TT_{\Delta_1,k_2,\Delta_2}).
\end{equation*}
\end{cor}

Actually, it follows from Lemma \ref{L3.4} and Remark
\ref{R3.5}that all sums in (ii), (iii) and (iv) are finite.

In this kind of situation, one can give a short interval version
of Theorem \ref{MainThm}. For each interval $I\subseteq [0,1]$ and
for each $\Delta=(\Delta_1,\dots,\Delta_h)\in (\N^*)^h$, let
\begin{equation*}
\begin{split}
& N_{Q,\odd}^I =\# \{ \gamma_0 <\dots <\gamma_h \
\mbox{consecutive in $\FF_{Q,\odd}$}\, ;\, \gamma_0 \in I\} \\ &
\qquad \qquad =
2\vert I\vert Q^2/\pi^2 +O(Q\log Q),  \\
& N_{Q,\odd}^I (\Delta) =\# \left\{
\begin{matrix} i ;\hspace{-6pt} & \gamma_i \in I,\
\gamma_i<\gamma_{i+1}<\dots <\gamma_{i+h} \
\mbox{consecutive in}\ \FF_{Q,\odd} \\
& \Delta(\gamma_{i+j-1},\gamma_{i+j})=\Delta_j,\ j=1,\dots,h
\end{matrix} \right\} .
\end{split}
\end{equation*}
Then the following result holds.

\medskip

\begin{thm}\label{Thm2}
Let $h\geq 1$, and assume that
$\Delta=(\Delta_1,\dots,\Delta_h)\in (\N^*)^h$ is such that only
finitely many non-vanishing terms appear on the right-hand side of
\eqref{I4}. Then for any interval $I\subseteq [0,1]$ we have
\begin{equation*}
\rho^I_{Q,\odd} (\Delta):=\frac{N_{Q,\odd}^I
(\Delta)}{N_{Q,\odd}^I} = \rho_\odd (\Delta)+O_{h,\varepsilon}
(Q^{-1/2+\varepsilon} )
\end{equation*}
for every $\varepsilon >0$.
\end{thm}

The main techniques of a proof involve the basic properties of
Farey fractions, the transformation $T$ from \eqref{I3}, and
estimates of Weil type for Kloosterman sums (see \cite{We},
\cite{Hoo}, \cite{Es}).

\bigskip

\section{Reduction of $N_{Q,\odd} (\Delta_1,\dots,\Delta_h)$}

We set throughout
\begin{equation*}
\Zf=\{ (a,b)\in \Z^2 \, ;\, \gcd (a,b)=1\},
\end{equation*}
and denote for any subset $\Omega$ of $\R^2$ and $Q\in \N^*$
\begin{equation*}
\begin{split}
\partial \Omega & =\mbox{\rm the boundary of $\Omega$},\qquad
Q \Omega =\{ (Qx,Qy)\, ;\, (x,y)\in \Omega \} ,\\
M(\Omega) & =\# (\Omega \cap \Z^2 ),\qquad M_\odd (\Omega)=\# \{
(x,y)\in \Omega \cap \Z^2 \, ;\, \mbox{\rm $x$ odd}\} ,\\ M_\even
(\Omega) & =\{ (x,y)\in \Omega \cap \Z^2 \, ;\, x\ \even \}
=M(\Omega)-
M_\odd (\Omega), \\
N(\Omega) & =\# (\Omega \cap \Zf ), \qquad N_\odd (\Omega)=\# \{
(x,y)\in \Omega \cap \Zf \, ;\,
\mbox{\rm $x$ odd}\}, \\
N_\even (\Omega) & =N(\Omega)-N_\odd (\Omega)=\{(x,y) \in \Omega
\cap \Zf \, ;\, \mbox{\rm $x$ even}\}, \\
N_{\odd,\odd} (\Omega) & =\{ (x,y)\in \Omega \cap \Zf \, ;\,
\mbox{\rm $x$ odd, $y$ odd}\},\\
N_{\odd,\even} (\Omega) & =\{ (x,y)\in \Omega \cap \Zf \, ;\,
\mbox{\rm $x$ odd, $y$ even}\},\\
N_{\even,\odd} (\Omega) & = \{ (x,y)\in \Omega \cap \Zf \, ;\,
\mbox{\rm $x$ even, $y$ odd} \} .
\end{split}
\end{equation*}

If $\gamma_{i_0}=\frac{a_{i_0}}{q_{i_0}} <\gamma_{i_0+1}=
\frac{a_{i_0+1}}{q_{i_0+1}} <\cdots <\gamma_{i_0+h}=
\frac{a_{i_0+h}}{q_{i_0+h}}$ are consecutive in $\FF_Q$, then (cf.
\cite{BCZ})
\begin{equation*}
\bigg( \frac{q_{i_0+r}}{Q} \, ,\frac{q_{i_0+r+1}}{Q} \bigg) = T^r
\bigg( \frac{q_{i_0}}{Q}\, ,\frac{q_{i_0+1}}{Q} \bigg) .
\end{equation*}

There is a one-to-one correspondence between $\Zf \cap Q
\TT_{k_1,\dots,k_r}$ and the set $\FF_{Q,k_1,\dots ,k_r}$ of
consecutive elements $\gamma_0 <\gamma_1 <\cdots <\gamma_r$ in
$\FF_Q$ with $\nu_Q (\gamma_{j-1})=k_j$, $j=1,\dots,r$, given by
\begin{equation*}
(q_0,q_1) \mapsto (\gamma_0,\gamma_1,\dots,\gamma_r),
\end{equation*}
where $(\gamma_0,\gamma_1)$ is the unique pair in $\FF_Q^<$ with
denominators $q_0$ and $q_1$, and $(\gamma_j,\gamma_{j+1})$ the
unique pair in $\FF_Q^<$ with denominators $Q T^j \big(
\frac{q_0}{Q},\frac{q_1}{Q} \big)$, $j=1,\dots,r$. This also shows
that the set
\begin{equation*}
\FF^{\odd,\odd/\even}_{Q,k_1,\dots,k_r} =\{
(\gamma_0,\dots,\gamma_r) \in \FF_{Q,k_1,\dots,k_r} \, ;\, q_0 \
\odd,\ q_1\ \odd/\even \}
\end{equation*}
has cardinality $N_{\odd,\odd/\even} (Q \TT_{k_1,\dots,k_r})$.

Suppose that $\gamma=\frac{a}{q} <\frac{a^{\prime
\prime}}{q^{\prime \prime}} =\gamma^{\prime \prime}$ are two
consecutive elements in $\FF_{Q,\odd}$, and that $\Delta (\gamma
,\gamma^{\prime \prime})=a^{\prime \prime} q-aq^{\prime
\prime}>1$. Since two fractions with even denominators cannot
occur as consecutive elements in $\FF_Q$, it follows that there is
precisely one fraction $\gamma^\prime =\frac{a^\prime}{q^\prime}$
in $\FF_Q$ such that $\gamma <\gamma^\prime <\gamma^{\prime
\prime}$ are consecutive in $\FF_Q$. One readily finds that (see
for example \cite[p.\,4]{Hay})
\begin{equation}\label{2.1}
\Delta (\gamma,\gamma^{\prime \prime})=\nu_Q (\gamma)= \bigg[
\frac{Q+q}{q^\prime} \bigg] =\kappa \bigg(
\frac{q}{Q},\frac{q^\prime}{Q} \bigg) .
\end{equation}

To summarize, suppose that $\gamma <\gamma^\prime <\gamma^{\prime
\prime} <\gamma^{\prime \prime \prime} < \gamma^{IV}$ are
consecutive in $\FF_Q$, and that $q$ is odd. Denote by
$q,q^\prime,\dots ,q^{IV}$, respectively, the denominator of
$\gamma ,\gamma^\prime,\dots ,\gamma^{IV}$. Denote also
\begin{equation*}
k_j=k_j (q,q^\prime)=\kappa \bigg( T^{j-1} \Big( \frac{q}{Q} ,
\frac{q^\prime}{Q} \Big)\bigg),  \qquad j\geq 1.
\end{equation*}
Then $q^{\prime \prime}=k_1 q^\prime -q$, $q^{\prime \prime
\prime}= k_2 q^{\prime \prime} -q^\prime$ etc. The following
situations may occur.
\begin{itemize}
\item[{\bf (O)}] $q^\prime$ is odd and thus $\Delta
(\gamma,\gamma^\prime)=1$. Next, it could be {\sl either} that
\begin{itemize}
\item[{\bf (OO)}] $q^{\prime \prime}$ is odd (if $k_1$ is even),
in which case $(\gamma^\prime,\gamma^{\prime \prime}) \in
\FF_{Q,\odd}^<$ and $\Delta (\gamma^\prime ,\gamma^{\prime
\prime})=1$, {\sl or} that \item[{\bf (OEO)}] $q^{\prime \prime}$
is even (if $k_1$ is odd), in which case $q^{\prime \prime \prime}
=k_2 q^{\prime \prime} -q^\prime$ is odd, $(\gamma^\prime ,
\gamma^{\prime \prime \prime})\in \FF_{Q,\odd}^<$, and $\Delta
(\gamma^\prime ,\gamma^{\prime \prime \prime})=k_2$.
\end{itemize}
\item[{\bf (E)}] $q^\prime$ is even, thus $q^{\prime \prime}$ is
odd, $(\gamma,\gamma^{\prime \prime})\in \FF^<_{Q,\odd}$ and
$\Delta (\gamma,\gamma^{\prime \prime})=k_1$. Next, we have {\sl
either} that
\begin{itemize}
\item[{\bf (EOO)}] $q^{\prime \prime \prime}$ is odd (if $k_2$ is
odd), in which case $(\gamma^{\prime \prime} ,\gamma^{\prime
\prime \prime})\in \FF_{Q,\odd}^<$ and $\Delta (\gamma^{\prime
\prime} ,\gamma^{\prime \prime \prime})=1$, {\sl or} that
\item[{\bf (EOEO)}] $q^{\prime \prime \prime}$ is even (if $k_2$
is even), in which case $q^{IV}=k_3 q^{\prime \prime \prime}
-q^{\prime \prime}$ will also be odd, $(\gamma^{\prime
\prime},\gamma^{IV})\in \FF_{Q,\odd}^<$ and $\Delta
(\gamma^{\prime \prime},\gamma^{IV})=k_3$.
\end{itemize}
\end{itemize}

This suggests that one may express $N_{Q,\odd}
(\Delta_1,\dots,\Delta_h)$ for any $h\geq 1$ by an inductive
procedure. Note first that
\begin{equation*}
N_{Q,\odd} (\Delta_1) =
\begin{cases}
\mbox{\small $\displaystyle \sum\limits_{k_1} N_{\odd,\odd}
(Q\TT_{k_1})+\sum\limits_{k_2} N_{\odd,\even} (Q\TT_{1,k_2})$}
 & \\ \mbox{\small $\displaystyle \qquad =N_{\odd,\odd} (Q\TT)+
N_{\odd,\even} (Q\TT_{\Delta_1})$}
& \mbox{\rm if $\Delta_1=1$} \\
\mbox{\small $\displaystyle \sum\limits_{k_2} N_{\odd,\even}
(Q\TT_{\Delta_1,k_2})=N_{\odd,\even} (Q\TT_{\Delta_1})$} &
\mbox{\rm if $\Delta_1 \geq 2$.}
\end{cases}
\end{equation*}

One may also express $N_{Q,\odd} (\Delta_1,\Delta_2)$ as
\begin{equation*}
\begin{cases}
\mbox{\small $\displaystyle \sum\limits_{k_1 \, \even}
N_{\odd,\odd} (Q\TT_{k_1})+\sum\limits_{k_1 \, \odd} N_{\odd,\odd}
(Q\TT_{k_1,1})$} & \\ \mbox{\small $\displaystyle \qquad
+\sum\limits_{k_2\, \odd} N_{\odd,\even} (Q\TT_{1,k_2})
+\sum\limits_{k_2 \, \even} N_{\odd,\even} (Q\TT_{1,k_2,1})$} &
\mbox{\rm if $\Delta_1= \Delta_2=1$}\\
\mbox{\small $\displaystyle \sum\limits_{k_1\,\odd} N_{\odd,\odd}
(Q\TT_{k_1,\Delta_2}) +\sum\limits_{k_2\, \even} N_{\odd,\even}
(Q\TT_{1,k_2,\Delta_2})$}
& \mbox{\rm if $\Delta_1=1,\Delta_2 \geq 2$} \\
\mbox{\small $\displaystyle \sum\limits_{k_2\, \odd}
N_{\odd,\even} (Q\TT_{\Delta_1,k_2})+\sum\limits_{k_2\, \even}
N_{\odd,\even} (Q\TT_{\Delta_1,k_2,1})$}
& \mbox{\rm if $\Delta_1 \geq 2,\Delta_2=1$} \\
\mbox{\small $\displaystyle \sum\limits_{k_2\, \even}
N_{\odd,\even} (Q\TT_{\Delta_1,k_2,\Delta_2})$} & \mbox{\rm if
$\Delta_1,\Delta_2 \geq 2$.}
\end{cases}
\end{equation*}

For $h\geq 2$, $\rho_{Q,\odd} (\Delta_1,\dots,\Delta_h)$ is
expressed in

\medskip

\begin{prop}\label{P2.1}
Assume that $h\geq 2$ and $\Delta =(\Delta_1,\dots,\Delta_h) \in
(\N^*)^h$. Then
\begin{equation}\label{2.2}
\rho_{Q,\odd} (\Delta) =\frac{1}{N_{Q,\odd}} \sum_{w\in
{\mathfrak{L}}_h \cap {\mathfrak{S}}_\Delta} N_{\odd,o(v_1)}
(Q\TT_{k_1,\dots,k_{\vert w\vert -1}}).
\end{equation}
\end{prop}

\bigskip

\section{Estimating $N_{\odd,\odd} (\Omega)$ and $N_{\odd,\even} (\Omega)$}

For a bounded region $\Omega$ in $\R^2$ with rectifiable boundary
and a function $f$ defined on $\Omega$, we set
\begin{equation*}
\begin{split}
& S_f (\Omega)=\sum\limits_{(a,b)\in \Omega \cap \Z^2}
\hspace{-8pt} f(a,b), \qquad S_f^\prime (\Omega)=
\sum\limits_{(a,b)\in \Omega \cap \Zf} \hspace{-8pt} f(a,b),\\
& S_{f,\odd / \even} (\Omega)= \sum_{\substack{(a,b)\in \Omega
\cap \Z^2 \\ a\, \odd / \even}} \hspace{-8pt} f(a,b), \qquad
S_{f,\odd/\even}^\prime = \sum_{\substack{(a,b)\in \Omega \cap \Zf
\\ a\, \odd / \even}}
\hspace{-8pt} f(a,b),\\
& S_{f,\odd,\odd/\even}^\prime (\Omega) =
\sum_{\substack{(a,b)\in \Omega \cap \Zf \\
a\, \odd,\,b\, \odd/\even}}
\hspace{-8pt} f(a,b) ,\\
& \| Df\|_{L^\infty (\Omega)} =\sup_{(x,y)\in \Omega} \bigg(
\bigg| \frac{\partial f}{\partial x} (x,y) \bigg| +\bigg|
\frac{\partial f}{\partial y} (x,y) \bigg| \bigg) .
\end{split}
\end{equation*}

\medskip

\begin{lem}\label{L3.1}
Let $R_1,R_2 >0$, and let $R\geq \min (R_1,R_2)$. Then for any
region $\Omega \subseteq [0,R_1 ]\times [0,R_2 ]$ and any function
$f$ which is $C^1$ on $\Omega$, we have
\begin{itemize}
\item[(i)] $\quad \displaystyle S_{f,\odd}^\prime
(\Omega)=\frac{4}{\pi^2} \iint\limits_\Omega f(x,y)\, dx\, dy
+O(A_{f,R,\Omega}).$ \item[(ii)] $\quad \displaystyle
S_{f,\odd,\odd/\even}^\prime (\Omega)= \frac{2}{\pi^2}
\iint\limits_\Omega f(x,y)\, dx\, dy +O(A_{f,R,\Omega}).$
\item[(iii)] $\quad \displaystyle S^\prime_{f,\even,\odd} (\Omega)
= \frac{2}{\pi^2} \iint\limits_\Omega f(x,y)\, dx\, dy+
O(A_{f,R,\Omega}),$
\end{itemize}
where
\begin{equation*}
\begin{split}
A_{f,R,\Omega} =\frac{\| f\|_{L^1 (\Omega)}}{R} & +\|
Df\|_{L^\infty (\Omega)} \area (\Omega) \log R \\ & + \|
f\|_{L^\infty (\Omega)} \big( R+\lh (\partial \Omega)\log R\big).
\end{split}
\end{equation*}
\end{lem}

{\sl Proof.} (i) It is well known (see e.g. \cite[Lemma\,1]{BCZ})
that
\begin{equation*}
S_f(\Omega)=\iint\limits_\Omega f(x,y)\, dx\, dy+ O(B_{f,\Omega}),
\end{equation*}
where
\begin{equation*}
B_{f,\Omega} =\| Df\|_{L^\infty (\Omega)} \area (\Omega) +\|
f\|_{L^\infty (\Omega)} \big( 1+\lh (\partial \Omega )\big).
\end{equation*}

Denoting $\Omega^\prime =\{ (x/2,y) \, ;\, (x,y)\in \Omega \}$, we
have that $S_{f,\even} (\Omega)$ - and eventually $S_{f,\odd}
(\Omega)$ can be expressed as
\begin{equation}\label{3.1}
\begin{split}
\sum\limits_{(a,b)\in \Omega^\prime \cap \Z^2} f(2a,b) &
=\iint\limits_{\Omega^\prime} f(2x,y)\, dx\, dy+O(
B_{f,\Omega^\prime}) \\ & =\frac{1}{2} \iint\limits_{\Omega}
f(x,y) \, dx\, dy+O( B_{f,\Omega}) .
\end{split}
\end{equation}

We now proceed to estimate $S^\prime_{f,\odd} (\Omega)$, which is
written as
\begin{equation}\label{3.2}
\begin{split}
\sum_{\substack{(a,b)\in \Omega \\ a\, \odd}} \hspace{-5pt} f(a,b)
& - \hspace{-8pt} \sum_{\substack{(a,b)\in \Omega /3 \\ a\, \odd}}
\hspace{-5pt} f(a,b) -\hspace{-8pt} \sum_{\substack{(a,b)\in
\Omega /5 \\ a\, \odd}} \hspace{-5pt} f(a,b) -\dots \\ &
=\sum_{\substack{1\leq n\leq R \\ n\, \odd}} \mu (n)
\sum_{\substack{(a,b)\in \Omega/n \\ a\, \odd}} f(na,nb) .
\end{split}
\end{equation}

The inner sum in \eqref{3.2}  is expressed by means of \eqref{3.1}
as
\begin{equation}\label{3.3}
\begin{split}
\frac{1}{2} & \iint\limits_{\Omega/n} f(nx,ny)\, dx\, dy \\ & +
O\Bigg( \frac{\| Df\|_{L^\infty (\Omega)} \area (\Omega)}{n} +\|
f\|_{L^\infty (\Omega)} \Big( 1+\frac{\lh (\partial \Omega)}{n}
\Big) \bigg) .
\end{split}
\end{equation}

Changing $(nx,ny)$ to $(x,y)$ in the double integral above and
summing over $n$, we infer from \eqref{3.2} and \eqref{3.3} that
\begin{equation*}
\begin{split}
S^\prime_{f,\odd} (\Omega) & =\frac{1}{2} \sum_{\substack{1\leq
n\leq R \\ n\, \odd}} \frac{\mu (n)}{n^2} \iint\limits_\Omega
f(x,y)\, dx\, dy +O\big( \|
Df\|_{L^\infty (\Omega)} \area (\Omega) \log R \big) \\
& +O\Big( \| f\|_{L^\infty (\Omega)} \big( R+\lh (\partial \Omega)
\log R \big) \Big).
\end{split}
\end{equation*}

The equality (i) now follows from
\begin{equation*}
\sum_{\substack{1\leq n\leq R \\ n\, \odd}} \frac{\mu (n)}{n^2} =
\frac{8}{\pi^2} +O\bigg( \frac{1}{R} \bigg) .
\end{equation*}

The equality (ii) follows by combining (i) with
\begin{equation*}
S^\prime_{f,\odd,\even} (\Omega) =
\sum_{\substack{(a,b)\in \Omega^{\prime \prime} \cap \Zf \\
a\, \odd}} \hspace{-8pt} f(a,2b),
\end{equation*}
where we set $\Omega^{\prime \prime} =\{ (x,y/2)\, ;\, (x,y)\in
\Omega \}$, and then using
\begin{equation*}
\iint\limits_{\Omega^{\prime \prime}} f(x,2y)\, dx\, dy=
\frac{1}{2} \iint\limits_\Omega f(x,y)\, dx\, dy.
\end{equation*}

The equality (iii) now follows from symmetry. \qed

We need the following improvement of Lemma 1 in \cite{Hay}.

\medskip

\begin{cor}\label{C3.2}
Let $R_1,R_2 >0$, and let $R\geq \min (R_1,R_2)$. Then for any
region $\Omega \subseteq [0,R_1 ]\times [0,R_2 ]$ with rectifiable
boundary, we have
\begin{itemize}
\item[(i)] $\displaystyle \quad N_\odd (\Omega)=4 \area
(\Omega)/\pi^2 +O(C_{R,\Omega})$, \item[(ii)] $\displaystyle \quad
N_\even (\Omega) =2\area (\Omega)/\pi^2 + O(C_{R,\Omega})$,
\item[(iii)] $\displaystyle \quad N_{\odd,\even} (\Omega)= 2\area
(\Omega)/\pi^2 +O(C_{R,\Omega})$, \item[(iv)] $\displaystyle \quad
N_{\odd,\odd} (\Omega)= 2\area (\Omega)/\pi^2 +O(C_{R,\Omega})$,
\item[(v)] $\displaystyle \quad N_{\even,\odd} (\Omega)= 2\area
(\Omega)/\pi^2 +O(C_{R,\Omega}),$
\end{itemize}
where
\begin{equation*}
C_{R,\Omega}=\area (\Omega)/R +R+\lh (\partial \Omega) \log R.
\end{equation*}
\end{cor}

The following lemma is contained in \cite{BCZ}. We enclose the
proof for the reader's convenience.

\medskip

\begin{lem}\label{L3.3}
For any integers $k_1,\dots ,k_r \geq 1$, the set
$\TT_{k_1,\dots,k_r}$ is a convex polygon.
\end{lem}

{\sl Proof.} If for $(x,y)\in \R^2$ we define $L_0(x,y)=x$,
$L_1(x,y)=y$, and $L_{i+1}(x,y)=k_i L_i(x,y)-L_{i-1}(x,y)$ for
$i\geq 1$, then $\TT_{k_1,\dots,k_r}$ is defined by the following
inequalities:
\begin{equation*}
\begin{cases}
1\geq L_0(x,y), L_1(x,y),\dots ,L_{r+1}(x,y)>0, \\
L_0(x,y)+L_1(x,y),L_1(x,y)+L_2(x,y),\dots,
L_r(x,y)+L_{r+1}(x,y)>1.
\end{cases}
\end{equation*}
Because $L_0,L_1,\dots,L_{r+1}$ are linear functions, the set
$\TT_{k_1,\dots,k_r}$ is the intersection of finitely many convex
polygons. \qed

\medskip

\begin{lem}\label{L3.4}
{\em (i)} Let $r\geq 1$. Then, for any $m\geq c_r=4r+2$, we have
that all sets $T^{-i} \TT_m$, $i=0,1,\dots,r$, are convex.
Moreover,
\begin{equation*}
T^{-1} \TT_m \subset \TT_1,\qquad \bigcup_{i=2}^r T^{-i} \TT_m
\subset \TT_2 ,
\end{equation*}
and, for all  $(x,y)\in \TT_m$ and $i\in \{ 1,2,\dots,r\}$,
\begin{equation*}
T^{-i} (x,y)=\big( x-iy,x-(i-1)y\big) .
\end{equation*}

{\em (ii)} For any $m\geq c_r$,
\begin{equation*}
T\TT_m \subset \TT_1 , \qquad \bigcup_{i=2}^r T^i \TT_m
 \subset \TT_2 ,
\end{equation*}
and, for all $(x,y)\in \TT_m$ and $i\in \{ 2,\dots ,r\}$,
\begin{equation*}
T^i(x,y)=\big( (m+2-i)y-x,(m+1-i)y-x\big) .
\end{equation*}

{\em (iii)} Let $j\in \{ 1,\dots ,r\}$. Then
\begin{equation*}
 \lh ( \partial T^{j-1} \TT_{k_1,\dots ,k_r} )
 \ll_r \frac{1}{k_j}
\end{equation*}
uniformly in $k_1,\dots,k_{j-1},k_{j+1},\dots,k_r$ as $k_j
\rightarrow \infty$.
\end{lem}

{\sl Proof.} (i) In the beginning we follow closely the proof of
Lemma 5 in \cite{BCZ}. The inverse of the transformation $T$ is
given by
\begin{equation}\label{3.4}
T^{-1} (x,y)=\bigg( \bigg[ \frac{1+y}{x} \bigg] x-y,x\bigg),
\qquad (x,y)\in \TT  .
\end{equation}
Since $0\leq 1-y<x$, we also have $\big[ \frac{1-y}{x} \big] =0$
and thus, for all $(x,y)\in \TT$,
\begin{equation}\label{3.5}
\kappa \big( T^{-1} (x,y)\big) =\left[ \frac{1+\big[
\frac{1+y}{x}\big] x-y}{x} \right] =\bigg[ \frac{1+y}{x} \bigg]
+\bigg[ \frac{1-y}{x} \bigg] = \bigg[ \frac{1+y}{x} \bigg] .
\end{equation}

Consider next a fixed element $(x,y)\in \TT_m$ with $m\geq c_r$.
Since $m\geq 5$, we have
\begin{equation*}
m\leq \frac{1+x}{y} <m+1 \qquad \mbox{\rm and} \qquad
x>\frac{m-1}{m+1} \, .
\end{equation*}
This leads to
\begin{equation*}
1<\frac{1+y}{x} \leq \frac{1+\frac{1+x}{m}}{x} =\frac{x+m+1}{mx}<
\frac{\frac{m-1}{m+1}+m+1}{m\cdot \frac{m-1}{m+1}}
=\frac{m+3}{m-1} \leq 2,
\end{equation*}
showing that $\kappa \big( T^{-1}(x,y)\big)=1$, and - using also
\eqref{3.4}- that
\begin{equation}\label{3.6}
T^{-1}(x,y)=(x-y,x)\in \TT_1 .
\end{equation}

Next, the inequality $m\geq c_r$ gives
\begin{equation}\label{3.7}
1+\frac{2(2r-1)}{m}\leq 1+\frac{m-4}{m} \leq 1+\frac{m-3}{m+1}
=\frac{2(m-1)}{m+1} \, .
\end{equation}
Thus the inequalities $x>\frac{m-1}{m+1}$ and $y\leq \frac{2}{m}$,
fulfilled by $(x,y)\in \TT_m$ (see \cite[Figure\,1]{BCZ}), imply
in conjunction with \eqref{3.7} that $2x>1+(2i-1)y$ for all
$(x,y)\in \TT_m$ and $i\in \{ 2,\dots,r\}$, or equivalently that
\begin{equation*}
\frac{1+y}{x-(i-1)y}<2 \, \qquad i\in \{ 2,\dots ,r\} .
\end{equation*}
At the same time, it is clear that $\frac{1+y}{x-(i-1)y}>1$, so
that
\begin{equation}\label{3.8}
\bigg[ \frac{1+x-(i-2)y}{x-(i-1)y} \bigg] = 1+\bigg[
\frac{1+y}{x-(i-1)y} \bigg] =2, \qquad i\in \{ 2,\dots,r\} .
\end{equation}

For $i=2$, equalities \eqref{3.5}, \eqref{3.7} and \eqref{3.8}
give
\begin{equation*}
\begin{split}
& \kappa \big( T^{-2} (x,y)\big) =\bigg[ \frac{1+x}{x-y} \bigg] =2, \\
& T^{-2} (x,y) =\big( 2(x-y)-x,x-y\big) =(x-2y,x-y),
\end{split}
\end{equation*}
thus, by \eqref{3.5} and by \eqref{3.8} with $i=3$ we have
\begin{equation*}
\begin{split}
& \kappa \big( T^{-3} (x,y)\big) =\bigg[ \frac{1+x-y}{x-2y}
\bigg]=2\qquad \mbox{\rm and} \\
& T^{-3}(x,y)=\big( 2(x-2y)-x+y,x-2y\big) =(x-3y,x-2y).
\end{split}
\end{equation*}

Arguing by induction, it follows at once that, for all $i\in \{
2,\dots,r\}$,
\begin{equation*}
\begin{split}
& \kappa \big( T^{-i} (x,y)\big) =\bigg[
\frac{1+x-(i-2)y}{x-(i-1)y}
\bigg] =2 \qquad \mbox{\rm and} \\
& T^{-i} (x,y)=\big( x-iy,x-(i-1)y\big) .
\end{split}
\end{equation*}

As a consequence, $T^{-i}\TT_m$ is the quadrangle with vertices at
$\big( 1-\frac{2i}{m},1-\frac{2(i-1)}{m} \big)$, $\big(
1-\frac{2i}{m+1} ,1-\frac{2(i-1)}{m+1}\big)$, $\big(
1-\frac{2(i+1)}{m+2} ,1-\frac{2i}{m+2}\big)$, and $\big(
1-\frac{2(i+1)}{m+1} ,1-\frac{2i}{m+1} \big)$. This quadrangle is
obviously contained in $\TT_2$.

(ii) Let $(x,y)\in \TT_m$. Then $T(x,y)=(y,my-x)$, and so
\begin{equation*}
\kappa \big( T(x,y)\big)=\bigg[ \frac{1+y}{my-x} \bigg] \geq 1.
\end{equation*}

Since $m\leq \frac{1+x}{y} <m+1$ and $y\leq \frac{2}{m} \leq
\frac{1}{3}$, it follows that $(2m-1)y\geq 1+(2m+2)y-2>1+2x$. This
leads to $\frac{1+y}{my-x} <2$, and so we obtain $\kappa \big(
T(x,y)\big)=1$. Therefore,
\begin{equation*}
T^2(x,y)=\big( my-x,(m-1)y-x\big) .
\end{equation*}

On the other hand, $y\leq \frac{1+x}{m}<\frac{1+x}{m-i}$; whence
\begin{equation}\label{3.9}
2\leq \bigg[ \frac{1+(m+2-i)y-x}{(m+1-i)y-x} \bigg] =1+ \bigg[
\frac{1+y}{(m+1-i)y-x} \bigg], \qquad i\geq 1 .
\end{equation}

The inequality $m\geq 4r+2$ leads to $m-2r\geq (2r+1)x$, which is
equivalent to $(m+1)(1+2x)\leq (2m+1-2r)(1+x)$. Since $1+x
<(m+1)y$, we infer that $1+2x<(2m+1-2r)y\leq (2m+1-2i)y$. That is,
\begin{equation}\label{3.10}
\frac{1+y}{(m+1-i)y-x} <2, \qquad i\in \{ 1,\dots ,r\} .
\end{equation}

By \eqref{3.9} and \eqref{3.10}, we gather that
\begin{equation}\label{3.11}
\bigg[ \frac{1+(m+2-i)y-x}{(m+1-i)y-x} \bigg] =2, \qquad i\in \{
2,\dots ,r\} .
\end{equation}

Now we infer inductively that $T^i (x,y)\in \TT_2$, and that
\begin{equation*}
T^i (x,y)=\big( (m+2-i)y-x,(m+1-i)y-x\big), \qquad i\in \{ 2,\dots
,r\} .
\end{equation*}

(iii) We use the fact that if $\Omega_1$ and $\Omega_2$ are convex
polygons with $\Omega_1 \subseteq \Omega_2$, then $\lh (\partial
\Omega_1)\leq \lh (\partial \Omega_2)$. For $k_j >c_r$ this
yields, in conjunction with (i) and (ii),
\begin{equation*}
\lh (\partial T^{j-1} \TT_{k_1,\dots,k_r}) \leq \lh
 (\partial \TT_{k_j}) \ll_r \frac{1}{k_j}
\end{equation*}
uniformly in $k_1,\dots,k_{j-1},k_{j+1},\dots ,k_r$. \qed

\begin{remark}\label{R3.5}
Suppose that $(x,y)\in \TT_m$ with $m\geq 3$. Then $\frac{1+x}{y}
<m+1$ and $y\leq \frac{2}{m}$; hence
$\frac{1+y}{my-x}<\frac{1+2/m}{1-y} \leq \frac{1+2/m}{1-2/m}
=\frac{m+2}{m-2}$ and thus
\begin{equation*}
\bigcup\limits_{m\geq 6} T\TT_m  \subset \TT_1 ,\quad T(\TT_4 \cup
\TT_5)\subset \TT_1 \cup \TT_2 \quad \mbox{\rm and} \quad T\TT_3
\subset \TT_1 \cup \TT_2 \cup \TT_3 \cup \TT_4 .
\end{equation*}
If $(x,y)\in \TT_2$, then $y>\frac{1+x}{3}\geq
\frac{x}{2}+\frac{1}{6} \geq \frac{1}{3}
+\frac{1}{6}=\frac{1}{2}$, and so
\begin{equation*}
T\TT_2 \subset \TT_1 \cup \TT_2 \cup \TT_3 .
\end{equation*}

On the other hand, if $(x,y)\in \TT_m$, $m\geq 2$, then it follows
by the proof of Lemma \ref{L3.4} (i), that $\kappa \big(
T^{-1}(x,y)\big)<\frac{m+3}{m-1}$. Therefore,
\begin{equation*}
\bigcup\limits_{m\geq 5} T\TT_m \subset \TT_1 , \quad T(\TT_3 \cup
\TT_4 )\subset \TT_1 \cup \TT_2 \quad \mbox{\rm and} \quad T\TT_2
\subset \TT_1 \cup \TT_2 \cup \TT_3 \cup \TT_4 .
\end{equation*}
\end{remark}

Owing to the presence of the term $R$ in $C_{R,\Omega}$, we need
one more fact, already noticed (in a different form) in
\cite{BGZ}.

\medskip

\begin{lem}\label{L3.6}
Let $k\in \N^*$ and let $\DD$ be a subset of $\TT$. Then the
following equalities hold.
\begin{itemize}
\item[{\rm (i)}] \ For $k$ even:
\begin{equation*}
\begin{split}
& N_{\odd,\even} \big( Q(\TT_k \cap \DD)\big) =
N_{\even,\odd}  \big( Q T(\TT_k \cap \DD )\big)  .\\
&  N_{\even,\odd} \big( Q(\TT_k \cap \DD)\big) =
N_{\odd,\even}  \big( Q T(\TT_k \cap \DD )\big)  . \\
&  N_{\odd,\odd} \big( Q(\TT_k \cap \DD)\big) = N_{\odd,\odd}
\big( Q T(\TT_k \cap \DD )\big) .
\end{split}
\end{equation*}
\item[{\rm (ii)}] \  For $k$ odd:
\begin{equation*}
\begin{split}
&  N_{\odd,\even} \big( Q(\TT_k \cap \DD)\big) =
N_{\even,\odd} \big( Q T(\TT_k \cap \DD )\big) .\\
&  N_{\even,\odd} \big( Q(\TT_k \cap \DD)\big) =
N_{\odd,\odd} \big( Q T(\TT_k \cap \DD )\big) . \\
&  N_{\odd,\odd} \big( Q(\TT_k \cap \DD)\big) = N_{\odd,\even}
\big( Q T(\TT_k \cap \DD )\big) .
\end{split}
\end{equation*}
\end{itemize}
\end{lem}

{\sl Proof.} We denote by $T_k$ the linear transformation defined
on $\R^2$ by $T_k (x,y)=(y,ky-x)$. Assume that $k$ is even and let
$(a,b)\in Q(\TT_k \cap \DD)$. Then $T\big( \frac{a}{Q},\frac{b}{Q}
\big)= \big( \frac{b}{Q},\frac{kb}{Q}-\frac{a}{Q} \big)$, so
\begin{equation*}
Q T(\TT_k \cap \DD)=\{ (b,kb-a) \, ;\, (a,b) \in Q\TT_k \cap Q\DD
\} =T_k \big( Q(\TT_k \cap \DD )\big) .
\end{equation*}
Moreover, since the matrix that defines $T_k$ is unimodular, the
elements of $\Zf \cap Q(\TT_k \cap \DD)$ are in $1-1$
correspondence with the elements of $\Zf \cap T_k \big( Q(\TT_k
\cap \DD )\big) =\Zf \cap Q\big( T(\TT_k \cap \DD )\big)$.
Besides, we see that $a$ is odd and $b$ is even if and only if $b$
is even and $kb-a$ is odd, implying that
\begin{equation*}
\begin{split}
\# \{ (a,b)\in \Zf & \cap Q(\TT_k \cap \DD ) \, ;\, a\, \odd,\ b\,
\even \} \\ & =\# \{ (c,d)\in \Zf \cap QT (\TT_k \cap \DD ) \, ;\,
c\, \even,\ d\, \odd \} .
\end{split}
\end{equation*}

The other five equalities follow in a similar way. \qed

\bigskip

{\sl Proof of Theorem {\em \ref{MainThm}.}} We wish to apply
Corollary \ref{C3.2} to $\Omega =Q\TT_{k_1,\dots,k_r}$. Note first
that, since $T$ is area-preserving, we have
\begin{equation*}
\area (\TT_{k_1,\dots,k_r}) \leq \area ( T^{-j+1} \TT_{k_j})
=\area (\TT_{k_j})\ll \frac{1}{k_j^3} \, \qquad j\in \{ 1,\dots
,r\} .
\end{equation*}
We claim that for every $j\in \{ 1,\dots ,r\}$, all the numbers
$N_{\odd,\odd} (Q\TT_{k_1,\dots,k_r})$, $N_{\odd,\even}
(Q\TT_{k_1,\dots,k_r})$, and $N_{\even,\odd}
(Q\TT_{k_1,\dots,k_r})$ can be expressed as
\begin{equation}\label{3.12}
\frac{2Q^2}{\pi^2} \, \area (\TT_{k_1,\dots,k_r})+ O_r \bigg(
\frac{Q}{k_j} \, \log Q\bigg)
\end{equation}
uniformly in $k_1,\dots,k_{j-1},k_{j+1},\dots,k_r$ as
$Q\rightarrow \infty$.

If $j\geq 2$, we apply Lemma \ref{L3.6} successively $j-1$ times:
to $k_1$ and $\DD=T^{-1} \TT_{k_2,\dots,k_r}$; to $k_2$ and
$\DD=T\TT_{k_1} \cap T^{-1} \TT_{k_3,\dots,k_r}$;$\dots$; and to
$k_{j-1}$ and $\DD=T^{j-2} \TT_{k_1,\dots,k_{j-2}} \cap T^{-1}
\TT_{k_j,\dots,k_r}$. This yields
\begin{equation*}
N_{\odd,\odd} (Q\TT_{k_1,\dots,k_r})=N_{\delta_1,\delta_2} (
QT^{j-1} \TT_{k_1,\dots,k_r})
\end{equation*}
for some pair $(\delta_1,\delta_2) \in \{
(\odd,\odd),(\odd,\even),(\even,\odd)\}$ that depends on
$k_1,\dots,k_{j-1}$. We may now apply Corollary \ref{C3.2} to
$\Omega =QT^{j-1} \TT_{k_1,\dots,k_r} \subseteq Q\TT_{k_j} \subset
[0,Q]\times \big[ 0,\frac{2Q}{k_j} \big]$, with $R\asymp
\frac{Q}{k_j}$, $\area (\Omega)\leq \area (Q\TT_{k_j})\ll
\frac{Q^2}{k_j^3}$, and (according to Lemma \ref{L3.4}) $\lh
(\partial \Omega) \ll_r \frac{Q}{k_j}$. Therefore, we gather that
$N_{\odd,\odd} (Q\TT_{k_1,\dots,k_r})$ is indeed given by
\eqref{3.12}. The same estimates are proved for $N_{\odd,\even}
(Q\TT_{k_1,\dots,k_r})$ and $N_{\even,\odd}
(Q\TT_{k_1,\dots,k_r})$ in a similar fashion.

We may now complete the proof of Theorem \ref{MainThm}. If $k_j
\geq c_r$, then we infer from Lemma \ref{L3.4} (i) that
$\TT_{k_1,\dots,k_r}=\emptyset$ unless $k_1=\dots =k_{j-2}=k_{j+2}
=\dots =k_r=2$ and $k_{j-1}=k_{j+1}=1$. On the other hand, we see
from \cite[Remark\,2.3]{BGZ} that $Q\TT_{k_1,\dots,k_r} \cap \Z^2
=\emptyset$ unless $\max (k_1,\dots,k_r) \leq 2Q$.

As a result, the only non-zero terms that may appear in the sum
from \eqref{2.2} arise from paths $w$ having all labels $k_j \leq
2Q$ and at most one $>c_{2h-1}$. Taking now also into account
\eqref{3.12}, the sum $\sum_{w\in {\mathfrak{L}}_h \cap
{\mathfrak{S}}_\Delta} \hspace{-3pt} N_{\odd,o(v_1)}
(Q\TT_{k_1,\dots ,k_{\vert w\vert -1}})$ can be expressed as
\begin{equation}\label{3.13}
\begin{split}
\frac{2Q^2}{\pi^2} \sum_{w\in {\mathfrak{L}}_h \cap
{\mathfrak{S}}_\Delta} & \hspace{-3pt} \area
(\TT_{k_1,\dots,k_{\vert w\vert -1}})+O_h
\left( \, \sum\limits_{k=1}^{2Q} \frac{Q\log Q}{k} \right) \\
& =\frac{2Q^2}{\pi^2} \sum_{w\in {\mathfrak{L}}_h \cap
{\mathfrak{S}}_\Delta} \hspace{-3pt} \area
(\TT_{k_1,\dots,k_{\vert w\vert -1}})+O_h (Q\log^2 Q) .
\end{split}
\end{equation}

The statement in Theorem \ref{MainThm} now follows from
Proposition \ref{P2.1}, \eqref{3.13}, and \eqref{1.1}. \qed

\bigskip

\section{Consecutive Farey fractions with odd denominators
in short intervals} For each interval $I\subseteq [0,1]$, and each
subset $\Omega \subseteq \R^2$, we set
\begin{equation*}
\Omega^I =\{ (a,b)\in \Omega \cap \Zf \, ;\, \bar{b}\in I_a \} ,
\end{equation*}
where $\bar{b}$ denotes the unique number in $\{ 1,\dots,a-1\}$
for which $b\bar{b}=1 \hspace{-3pt} \pmod{a}$. If
$I=[\alpha,\beta]$, then we also set $I_a
=[a(1-\beta),a(1-\alpha)]$.

For any function $f$ defined on $\Omega$, denote
\begin{equation*}
S^I_{f,\odd,\odd/\even} (\Omega)= \sum_{\substack{(a,b)\in \Omega
\cap \Zf \\ a\, \odd,\, b\, \odd/\even \\ \bar{b} \in I_a}}
\hspace{-8pt}  f(a,b).
\end{equation*}

The following analog of Proposition \ref{P2.1} holds and is
similarly proved.

\medskip

\begin{prop}\label{P4.1}
Let $h\geq 1$, and let $\Delta=(\Delta_1,\dots,\Delta_h) \in
(\N^*)^h$. Then, for any interval $I\subseteq [0,1]$,
\begin{equation*}
N^I_{Q,\odd} (\Delta)= \sum_{\substack{w\in {\mathfrak{L}}_h \cap
{\mathfrak{S}}_\Delta}} N_{\odd,o(v_1)} \big(
(Q\TT_{k_1,\dots,k_{\vert w\vert -1}})^I \big).
\end{equation*}
\end{prop}

\medskip

\begin{prop}\label{P4.2}
Assume that $\Omega \subseteq [0,R_1]\times [0,R_2 ]$ is a convex
region, and that $f$ is a $C^1$ function on $\Omega$. Then
$S^I_{f,\odd,\odd/\even} (\Omega)$ is given by
\begin{equation*}
\vert I\vert S^\prime_{f,\odd,\odd/\even} (\Omega) +O_\varepsilon
\big( \| f\|_{L^\infty (\Omega)} (R_2 \log R_1 +m_f
R_1^{1/2+\varepsilon} (R_1+R_2) )\big)
\end{equation*}
for every $\varepsilon >0$, where $m_f$ is an upper limit for the
number of intervals of monotonicity of the functions $y\mapsto
f(x,y)$.
\end{prop}

{\sl Proof.} The proof is similar to that of Lemma 8 in
\cite{BCZ}. As in \cite[(65)]{BCZ}, we write
\begin{equation}\label{4.1}
S_{f,\odd,\odd/\even}^I (\Omega)=S_1+S_2,
\end{equation}
where
\begin{equation}\label{4.2}
\begin{split}
S_1 & =\sum_{\substack{(a,b)\in \Omega \cap \Zf \\ a\, \odd,\, b\,
\odd/\even}} \hspace{-8pt} f(a,b)\sum\limits_{x\in I_a}
\frac{1}{a} \\ & = \sum_{\substack{(a,b)\in \Omega \cap \Zf \\ a\,
\odd ,\, b\, \odd/\even}} \hspace{-8pt} f(a,b) \, \frac{1}{a} \,
\big( \vert I_a \vert +O(1)\big) \\
& =\vert I\vert S^\prime_{f,\odd,\odd/\even} (\Omega) +O(\|
f\|_{L^\infty (\Omega)} R_2 \log R_1 )
\end{split}
\end{equation}
and
\begin{equation*}
S_2 =\sum_{\substack{(a,b)\in \Omega \cap \Zf \\ a\, \odd ,\, b\,
\odd/\even}} \hspace{-8pt} f(a,b) \sum\limits_{x\in I_a}
\frac{1}{a} \sum\limits_{l=1}^{a-1} e\bigg( \frac{l(\bar{b}-x)}{a}
\bigg) .
\end{equation*}

As in \cite[(67)]{BCZ}, we write
\begin{equation}\label{4.3}
S_2 =\sum_{\substack{a\in {\mathrm{pr}}_1 (\Omega) \\ a\, \odd}}
\frac{1}{a} \sum\limits_{l=1}^{a-1} \left( \, \sum\limits_{x\in
I_a} e\bigg( -\frac{lx}{a} \bigg) \right)
S_{f,\odd/\even,I_a^\prime} (l,a),
\end{equation}
where $I_a^\prime =\{ b\, ;\, (a,b)\in \Omega\}$ is an interval
for every $a$ in the projection ${\mathrm{pr}}_1 (\Omega)$ of
$\Omega$ on the first coordinate. Here, for any interval $J$, we
denote
\begin{equation}\label{4.4}
S_{f,\odd/\even,J} (l,a)=\sum_{\substack{b\in J \\
b\, \odd/\even \\ \gcd (a,b)=1}} \hspace{-8pt} f(a,b)\ e\bigg(
\frac{l\bar{b}}{a} \bigg)
\end{equation}
and
\begin{equation*}
S_{f,J}(l,a)=\sum_{\substack{b\in J \\ \gcd (a,b)=1}}
\hspace{-8pt} f(a,b) \ e\bigg( \frac{l\bar{b}}{a} \bigg) .
\end{equation*}

By \cite[Lemma\,9]{BCZ} we have
\begin{equation}\label{4.5}
\vert S_{f,J} (l,a)\vert \ll_\varepsilon R_{\Omega
,f,J,l,a,\varepsilon} ,
\end{equation}
where
\begin{equation*}
R_{\Omega,f,J,l,a,\varepsilon} =m_f \| f\|_{L^\infty (\Omega)}
\big( \vert J\vert a^{-1/2+\varepsilon} + a^{1/2+\varepsilon}
\big) \gcd (l,a)^{1/2} .
\end{equation*}

Writing now
\begin{equation*}
S_{f,\even,J} (l,a)=\sum_{\substack{c\in J/2 \\
\gcd (a,c)=1}} \hspace{-8pt} f(a,2c)\, e\bigg( \frac{\bar{2}\,
l\,\bar{c}}{a} \bigg) =S_{f_2,J/2} (\bar{2} l,a),
\end{equation*}
where $f_2(x,y)=f(x,2y)$, and then using \eqref{4.5} and
$S_{f,\odd,J} (l,a)=S_{f,J}(l,a)-S_{f,\even,J} (l,a)$, we infer
that
\begin{equation}\label{4.6}
\max \big( \vert S_{f,\even,J}(l,a)\vert,\vert S_{f,\odd ,J}
(l,a)\vert \big) \ll_\varepsilon R_{\Omega,f,J,l,a,\varepsilon} .
\end{equation}

As in \cite[(67)-(69)]{BCZ}, we infer - from \eqref{4.3},
\eqref{4.4}, \eqref{4.6}, from the fact that the inner sum in
\eqref{4.3} is a geometric progression which is $\ll \big(
\frac{a}{l},\frac{a}{a-l} \big)$, and from $\vert I_a^\prime \vert
\leq R_2$ - that
\begin{equation}\label{4.7}
\begin{split}
\vert S_2 \vert & \ll \sum\limits_{a=1}^{R_1} \, \frac{1}{a} \,
\sum\limits_{l=1}^{a-1} \, \frac{a}{l} \, \vert
S_{f,\odd/\even,I_a^\prime} (l,a)\vert \\ & \ll_\varepsilon m_f \|
f\|_{L^\infty (\Omega)} R_1^{1/2+\varepsilon} (R_1+R_2) .
\end{split}
\end{equation}

The desired conclusion follows now from \eqref{4.1}, \eqref{4.2}
and \eqref{4.7}. \qed

\medskip

\begin{cor}\label{C4.3}
\begin{equation*}
N_{\odd,\odd/\even} \big( (Q\TT_{k_1,\dots,k_r})^I \big)= \vert
I\vert N_{\odd,\odd/\even} (Q\TT_{k_1,\dots,k_r}) +O_\varepsilon
(Q^{3/2+\varepsilon}) .
\end{equation*}
\end{cor}

\medskip

Theorem \ref{Thm2} is now a consequence of Proposition \ref{P4.1}
and Corollary \ref{C4.3}.

\bigskip

{\bf Acknowledgments.} We are grateful to the referee for careful
reading of the manuscript and pertinent suggestions that led to
the improvement of this paper.

\bigskip

\bibliographystyle{amsplain}

\begin{thebibliography}{4}
\bibitem{ABCZ} V. Augustin, F.P. Boca, C. Cobeli, A. Zaharescu,
\emph{The $h$-spacing distribution between Farey points}, Math.
Proc. Camb. Phil. Soc. {\bf 131} (2001), 23--38.

\bibitem{BCZ0} F.P. Boca, C. Cobeli, A. Zaharescu, \emph{Distribution of
lattice points visible from the origin}, Comm. Math. Phys. {\bf
213} (2000), 433--470.

\bibitem{BCZ} F.P. Boca, C. Cobeli and A. Zaharescu, \emph{A conjecture
of R.R. Hall on Farey points}, J. Reine Angew. Math. {\bf 535}
(2001), 207--236.

\bibitem{BGZ} F.P. Boca, R.N. Gologan, A. Zaharescu, \emph{On the index of
Farey sequences}, preprint arXiv math.NT/0201044, Quart. J. Math.
Oxford Ser. (2) {\bf 53} (2002), 377--391.

\bibitem{Es}
T. Esterman, \emph{On Kloosterman's sums}, Mathematika {\bf 8}
(1961), 83--86.

\bibitem{Fra} J. Franel, \emph{Les suites de Farey et le probl\` eme de
nombres premiers}, G\" ottinger Nachr., 1924, 198--201.

\bibitem{Hall1} R.R. Hall, \emph{A note on Farey series},
J. London Math. Soc. {\bf 2} (1970), 139--148.

\bibitem{Hall2} R.R. Hall, \emph{On consecutive Farey arcs II},
Acta Arith. {\bf 66} (1994), 1--9.

\bibitem{HT} R.R. Hall, G. Tenenbaum, \emph{On consecutive Farey arcs},
Acta Arith. 44 (1984), 397--405.

\bibitem{HSZ} R.R. Hall, P. Shiu, \emph{The index of a Farey sequence},
Michigan Math. J. {\bf 51} (2003), 209--223.

\bibitem{Hay} A. Haynes, \emph{A note on Farey fractions with odd
denominators}, J. Number Theory {\bf 98} (2003), 89--104.

\bibitem{Hoo}
C. Hooley, \emph{An asymptotic formula in the theory of numbers},
Proc. London Math. Soc. {\bf 7} (1957), 396--413.

\bibitem{Hux} M.N. Huxley, \emph{The distribution of Farey points I},
Acta Arith. {\bf 18} (1971), 281--287.

\bibitem{HZ} M.N. Huxley, A. Zhigljavsky,
\emph{On the distribution of Farey fractions and hyperbolic
lattice points}, Period. Math. Hungarica {\bf 42} (2001),
191--198.

\bibitem{La} E. Landau, \emph{Bemerkungen zu der vorstehenden
Abhandlung von Herrn Franel}, G\" ottinger Nachr. 1924, 202--206.

\bibitem{We}
A. Weil A, \emph{On some exponential sums}, Proc. Nat. Acad. Sci.
U.S.A. 34 (1948), 204--207.

\end{thebibliography}

\end{document}